\documentclass[12pt]{article}

\usepackage[latin1]{inputenc}
\usepackage{amsmath,amsthm,amssymb,amscd,}
\newtheorem{thm}{Theorem}[section]
 \newtheorem{cor}[thm]{Corollary}
 \newtheorem{lem}[thm]{Lemma}
 \newtheorem{prop}[thm]{Proposition}
 \theoremstyle{definition}
 \newtheorem{df}[thm]{Definition}
 \theoremstyle{remark}
 \newtheorem{rem}[thm]{Remark}
 
 \numberwithin{equation}{section}

\def\be#1 {\begin{equation} \label{#1}}
\newcommand{\ee}{\end{equation}}
\def\dem {\noindent {\bf Proof : }}

\def\sqw{\hbox{\rlap{\leavevmode\raise.3ex\hbox{$\sqcap$}}$%
\sqcup$}}
\def\findem{\ifmmode\sqw\else{\ifhmode\unskip\fi\nobreak\hfil
\penalty50\hskip1em\null\nobreak\hfil\sqw
\parfillskip=0pt\finalhyphendemerits=0\endgraf}\fi}

\newcommand{\mb}{\medskip\noindent}
\newcommand{\gb}{\bigskip\noindent}
\newcommand{\R}{\mathbb R}

\newcommand{\Z}{\mathbb Z}

\newcommand{\s}{\mathcal S}

\newcommand{\BS}{\overline{BS}}

\title{A bilinear pseudodifferential calculus.}
\author{Fr\'ed\'eric Bernicot \\
 \small Universit\'e de Paris-Sud, Orsay et CNRS 8628, 91405 Orsay Cedex, France \\
\small {\em E-mail address:} {Frederic.Bernicot@math.u-psud.fr} }

\begin{document}
\maketitle

\begin{abstract}
 In this paper, we are interested in the construction of a bilinear pseudodifferential calculus. We define some symbolic classes which contains those of Coifman-Meyer. These new classes allow us to consider operators closely related to the bilinear Hilbert transform. We give a description of the action of our bilinear operators on Sobolev spaces. These classes also have a ``nice'' behavior through the transposition and the composition operations that we will present. 

\mb {\bf Key words} : multilinear pseudodifferential calculus, time-frequency analysis, asymptotic expansion. 

\mb {\bf MSC classification} : 47G30-42B15-42C10.
\end{abstract}

\tableofcontents

\section{Introduction}

\subsection{The general approach of a bilinear pseudodifferential calculus.}

In the linear theory, the pseudodifferential calculus brings powerful tools to study Partial Differential Equations.
Its construction began with the study of the "classical" symbolic classes $(S^m_{1,0})_{m\in \R}$. Theses classes produce linear operators acting on Sobolev spaces. To use these operators, it is useful to understand a few functional rules, mainly the behavior of the adjointness and the composition operations on these symbolic classes. These two operations are described in the linear pseudodifferential calculus by asymptotic formulas (see \cite{AG}). Aiming to generalize, many people have searched to get the largest symbolic classes, that retain these properties. Certain works treat exotic classes of linear symbols (for
example see the book of R. Coifman and Y. Meyer \cite{cm}). \\
Nowadays some multilinear estimates appear to study the nonlinear terms in some P.D.E. (for example in the
``I-method''). This motivate us to define and to build a multilinear pseudodifferential calculus. In addition, we have shown in \cite{B} continuity in Sobolev spaces for new singular bilinear operators. That is why in this paper, we are interested to continue the definition and the study of a bilinear pseudodifferential calculus, which were started in
\cite{bipseudo,BT} by A. Bényi, A. Nahmod and R. Torres. Let us first describe the already known results about bilinear pseudodifferential operators. \\
We know that to a bilinear operator $T$, which is continuous
from $\s(\R) \times \s(\R)$ into $\s'(\R)$, we can associate a symbol $\sigma\in\s(\R^3)'$ such that~:
\begin{align}
 \forall f,g\in\s(\R),\qquad T(f,g)(x) & = \int e^{ix(\alpha+\beta)} \sigma(x,\alpha,\beta)\widehat{f}(\alpha) \widehat{g}(\beta) d\alpha d\beta \label{operator} \\
  & :=T_\sigma(f,g)(x). \nonumber
\end{align}
Our goal is to obtain the largest class of bilinear symbols $\sigma$ such that the operators $T_\sigma$ verify
continuities in Lebesgue and Sobolev spaces. In addition, we would like to understand the behavior of the adjointness and composition operations on these operators. 


\gb Let us recall the first classes of symbols $BS^m_{\rho,\delta}$. 

\begin{df} For all real $m \in \R$ and reals $0\leq \rho,\delta\leq 1$, a function $\sigma\in C^\infty(\R^3)$ belongs to the class $BS^m_{\rho,\delta}$ if~:
\be{BS1} \forall a,b,c\geq 0,\quad \left|  \partial_x^a \partial_\alpha^b \partial_\beta^c \sigma(x,\alpha,\beta) \right| \lesssim \left(1+|\alpha|+|\beta|\right)^{m+\delta a-\rho(b+c)}. \ee  
\end{df}

\mb The class $BS^0_{1,0}$ goes back to the work of R. Coifman and Y. Meyer \cite{cm} while for $m\neq 0$
 the classes $BS^m_{1,0}$  were defined by A. Benyi and R. Torres in \cite{BT} and started to be studied 
 by A. Benyi, A. Nahmod and R. Torres in \cite{bipseudo}.
 
\mb The boundedness properties of the class $BS^0_{1,0}$ are nowadays well understood and are given by the following result~:
 
\begin{thm} Let $0<p,q,r <\infty$ be exponents such that
$$ \frac{1}{r}=\frac{1}{p}+\frac{1}{q} \qquad \textrm{and} \qquad 1<p,q<\infty.$$
Then for all symbol $\sigma \in BS^0_{1,0}$, the operator $T_\sigma$ is continuous from $L^p(\R) \times L^q(\R)$ into $L^r(\R)$.
\end{thm}
 
\mb The result was proved by R. Coifman and Y. Meyer for $r>1$ in \cite{cm}, while for $r>1/2$  it was obtained 
 by C. Kenig and E. Stein in \cite{KS} and L. Grafakos and R. Torres in \cite{GT} via interpolation and a weak-type 
 end-point estimate for $r=1/2$.

\mb In addition A. Bényi and R.Torres have shown in \cite{BT} that the class $BS^{0}_{1,0}$ is closed under transposition. A. Bényi have obtained in \cite{Benyi} similar results for the  ``more exotic'' classes $BS^0_{1,\delta}$ with $0\leq \delta<1$. In fact, all these symbols generate bilinear Calder\'on-Zygmund operators and also their study are closely related to the multilinear Calder\'on-Zygmund theory.

\gb Nowadays, many people are interested by far more singular operators, which are completely outside the Calder\'on-Zygmund theory. The prototype of these new operators is the bilinear Hilbert transform, which appeared in the work of A. Calder\'on (\cite{calderon1,calderon2} in 60-70's). His famous conjecture, about the continuity of this one from $L^2(\R) \times L^\infty(\R)$ into $L^2(\R)$, was later solved by M. Lacey and C. Thiele (\cite{LT2,LT1,LT3,LT4} in 90's). Then C. Muscalu, T. Tao and C. Thiele (in \cite{MTT2}) and independently J. Gilbert and A. Nahmod (in \cite{GN}) have extended the proof to study a new class of bilinear operators. We are also interested to contruct a bilinear pseudodifferential calculus, with symbolic classes containing this kind of operators. Such a bilinear calculus was already defined by A. B\'enyi, A. Nahmod and R. Torres in \cite{bipseudo}, so we first recall their definitions.

\begin{df} For all real $m \in \R$, all reals $0\leq \rho,\delta\leq 1$ and all angle $\theta \in ]-\pi/2,\pi/2[$, a function $\sigma\in C^\infty(\R^3)$ belongs to the class $BS^m_{\rho,\delta;\theta}$ if~:
\be{BS2} \forall a,b,c\geq 0,\quad  \left|  \partial_x^a \partial_\alpha^b \partial_\beta^c \sigma(x,\alpha,\beta) \right| \lesssim \left(1+|\beta-\tan(\theta)\alpha|\right)^{m+\delta a-\rho(b+c)}. \ee 
\end{df}

\mb Also we have replaced in (\ref{BS1}) the quantity $|\alpha|+|\beta|=d((\alpha,\beta),0)$ by the lower quantity $\beta-\tan(\theta)\alpha=d((\alpha,\beta),\Delta)$, where $\Delta$ is the line 
$\Delta:=\{(\alpha,\beta),\ \beta=\tan(\theta)\alpha\}$ in the frequency plane. It is obvious that for every angle $\theta$, $BS_{1,0}^{0} \subset BS_{1,0;\theta}^{0}$. \\
The previously mentioned papers \cite{MTT2,GN} deal with the case where the symbol $\sigma$ belongs to the class $BS^0_{1,0;\theta}$ and is $x$-independent. Their main result is the following one~:
\begin{thm} \label{aze} Let $\Delta$ be an nondegenerate line of the frequency plane. This means that
$$ \Delta:=\left\{(\alpha,\beta), \beta-\tan(\theta) \alpha=0 \right\} $$
with $\theta \in ]-\pi/2,\pi/2[ \setminus\{-\pi/4,0\}$. \\
Let $p,q$ be exponents such that
$$1< p,q \leq \infty \quad \textrm{and} \quad  0<\frac{1}{r}=\frac{1}{q}+\frac{1}{p} <\frac{3}{2}.$$
Then for all $x$-independent symbol $\sigma\in BS^0_{1,0;\theta}$, the operator $T_\sigma$
is continuous from $L^p(\R) \times L^q(\R)$ in $L^r(\R)$.
\end{thm}
\mb In \cite{t1bilineaire}, the authors have extended this result in a particular case : if the symbol
$\sigma$ verifies $\sigma(x,\alpha,\beta)=\tau(x,\beta-\tan(\theta)\alpha)$ with a symbol $\tau$ satisfying for all
$a,b\geq 0$ \be{cond4} \left|\partial^{a}_x \partial^b_\xi \tau(x,\xi)\right| \lesssim \left(1+|\xi|\right)^{-b},  \ee
then the bilinear operator $T_\sigma$ satisfies the same continuities. 
All these symbols could be much more singular than those of the class $BS^0_{1,0}$. The arguments, used in the proof of these continuities, are based on a very nice and sharp time-frequency analysis. \\
In  \cite{bipseudo}, the authors ask a question : Is it possible to obtain a result like Theorem  \ref{aze} for $x$-dependent symbols ? In addition they particularly study the special class $BS^{0}_{1,0;-\pi/4}$ and they study the action of the duality on the operators associated to the class $BS^{0}_{1,0;\theta}$.

\gb In this paper,  we want to continue the construction of such a bilinear pseudodifferential calculus. In \cite{B}, we have positively answered to the previous question. So the operators associated to the classes $BS^0_{1,0;\theta}$ act on Sobolev spaces. However theses symbolic classes are not invariant by composition with linear pseudodifferential operators (as explained in Theorem \ref{thm:composition}). We want also to consider larger symbolic classes closely related to these classes $BS^m_{1,0;\theta}$, and we would like to decompose the order $m$
in two orders for the two frequency variables $\alpha$ and $\beta$. So we will construct some larger classes which verify continuities on Sobolev spaces and some functional invariances. 

\subsection{Notations and our main results.}

For notation, we denote the norm in $L^p(E)$ for any measurable set $E\subset \R$ by $\|\ \|_{p,dx,E}$ (or $\|\
\|_{p,E}$ if there is no confusion for the variable). 

\mb In this subsection, we begin to define our new classes of bilinear symbols. First we recall the
classical linear classes.

\begin{df}  For $m\in \R$, we denote $S^{m}_{1,0}$ the classical set of linear pseudodifferential symbol of order $m$.
$$ S^{m}_{1,0}:=\left\{ \sigma\in C^\infty(\R^2),\ \left| \partial_x^a \partial_\alpha ^b \sigma(x,\alpha) \right| \lesssim (1+|\alpha|)^{m-b}\ \right\}.$$
\end{df}

\mb We use the same idea to define the bilinear classes~:

\begin{df} For $m_1,m_2$ two reals, we define the class of bilinear pseudodifferential symbols of order $(m_1,m_2)$. Let $\theta$ be an angle and $\sigma$ be a $C^\infty(\R^3)$ function. \\
We set that $\sigma$ belongs to the class $\BS_{1,0;\theta}^{m_1,m_2}$ if and only if for all $a,b,c\geq 0$
\begin{align*}
\lefteqn{\left| \partial_x^a \partial_\alpha ^b \partial_\beta^c \sigma(x,\alpha,\beta) \right| \lesssim \left(1+|\alpha| \right)^{m_1}\left(1+|\beta| \right)^{m_2} } & & \\
& & \left(1+\min\{|\alpha|,|\beta-\tan(\theta)\alpha| \}\right)^{-b} \left(1+\min\{|\beta|,|\beta-\tan(\theta)\alpha|\}\right)^{-c}.
\end{align*}
We set that $\sigma$ belongs to the class $\BS^{m_1,m_2,1}_{1,0;\theta}$ if and only if for all $a,b,c\geq 0$ 
\begin{align*}
\lefteqn{ \left| \partial_x^a \partial_\alpha ^b \partial_{\beta-\alpha}^c \sigma(x,\alpha,\beta) \right| \lesssim \left(1+|\alpha+\beta| \right)^{m_1}\left(1+|\beta| \right)^{m_2} } & & \\
 & &   \left(1+\min\{|\alpha+\beta|,|\beta-\tan(\theta)\alpha| \}\right)^{-b} \left(1+\min\{|\beta|,|\beta-\tan(\theta)\alpha|\}\right)^{-c}.
 \end{align*}
Finally we set that $\sigma$ belongs to the class $\BS^{m_1,m_2,2}_{1,0;\theta}$ if and only if for all $a,b,c\geq 0$ 
\begin{align*}
\lefteqn{ \left| \partial_x^a \partial_{\alpha-\beta} ^b \partial_\beta^c \sigma(x,\alpha,\beta) \right| \lesssim \left(1+|\alpha| \right)^{m_1}\left(1+|\alpha+\beta| \right)^{m_2} } & & \\
 & & \left(1+\min\{|\alpha|,|\beta-\tan(\theta)\alpha| \}\right)^{-b} \left(1+\min\{|\alpha+\beta|,|\beta-\tan(\theta) \alpha|\}\right)^{-c}.
 \end{align*}
\end{df}

\mb In the three previous definitions, the most important term is the additional decay obtained when we differentiate the symbol with respect to the frequency variables. With the notations of \cite{bipseudo}, our class $\BS_{1,0;\theta}^{m,m}$ contains the class
$BS^{m}_{1,0;\theta}$.

\mb
\begin{rem} For $\theta \in ]-\pi/2,\pi/2[-\{0\}$ we have the following equivalence
$$ \sigma(x,\alpha,\beta)\in \BS^{m_1,m_2}_{1,0,\pi/4} \Longleftrightarrow \sigma(x,\tan(\theta)\alpha,\beta) \in \BS_{1,0;\theta}^{m_1,m_2}.$$ We have the same equivalence for the other classes $\BS^{m_1,m_2,1}_{1,0;\theta}$ and $\BS^{m_1,m_2,2}_{1,0;\theta}$.
\end{rem}

\gb
\begin{rem} In the definition of our symbolic classes, the term $\min\{|\alpha|,|\beta-\tan(\theta)\alpha|\}$ corresponds to the distance in the frequency plane $d((\alpha,\beta),C_{\theta})$ between the cone $C_\theta$ (composed of the two lines
$\beta=\tan(\theta)\alpha$ and $\alpha=0$) and the point $(\alpha,\beta)$. In \cite{GN,GN2} J.Gilbert and A.Nahmod have studied the case where the cone is
nondegenerate (i.e. the two lines of the cone are nondegenerate). Here the cone $C_\theta$ has a degenerate line. However the difference is that in our definition of
the class $\BS^{m_1,m_2}_{1,0;\theta}$, the cone is not the same for the differentiations with respect to $\alpha$ and to $\beta$. We have two
different degenerate cones for the two frequency variables. This will allow us to obtain continuities in Lebesgue and Sobolev spaces for the associated operators.
\end{rem}

\begin{df} For $m$ a real, we set $m_+=(|m|+m)/2$. We denote $L^{p}=L^{p}(\R)$ for the Lebesgue spaces. The Sobolev space $W^{m,p}=W^{m,p}(\R)$ is defined as the set of distributions $f\in \s'(\R)$ such that
$$ \|f\|_{W^{n,p}}:=\|J_m(f)\|_{L^p},$$
where $J_m:=\left(Id-\Delta \right)^{m/2}$.
\end{df}

\mb Now we come to our main results. With the same ideas that we used in \cite{B}, we will first study the action of our pseudodifferential operators on Sobolev spaces and we will prove the following result~:

\begin{thm} \label{thm:central} Let $\sigma$ be a symbol of order $(m_1,m_2)$ (in one of the three classes $\BS^{m_1,m_2}_{1,0;\theta}$, $\BS^{m_1,m_2,1}_{1,0;\theta}$ or $\BS^{m_1,m_2,2}_{1,0;\theta}$ for $ \theta\in]-\pi/2,\pi/2[ \setminus \{0,-\pi/4\}$). Let $p,q,r$ be exponents satisfying
$$ 0<\frac{1}{r}=\frac{1}{p}+\frac{1}{q}<\frac{3}{2} \quad \textrm{and} \quad 1<p,q\leq \infty.$$
Then for all real $s\geq 0$ and $\epsilon>0$, the operator $T_\sigma$ is continuous from $W^{s+\epsilon+(m_1)_+,p}\times W^{s+\epsilon+(m_2)_+,q}$ into $W^{s,r}$.
\end{thm}

\mb In the previous Theorem, the $\epsilon$ is necessary, due to the assumed Marcinkiewicz's conditions on the symbol (see \cite{GK}). We will give other continuities, with no loss of regularity, by using the Modulation spaces and we will give some extra condition to allow us to get the previous continuities with $\epsilon=0$. 

 \mb After this main result,  we will describe in Section \ref{rules} two rules of the bilinear pseudodifferential calculus~:
 
\begin{thm} \label{thm:dualite} For $T$ a bilinear operator, we write $T^{*1}$ and $T^{*2}$ for its two adjoints ( with respect to the $\alpha$ and to the $\beta$ variable). Let $m_1,m_2\in\R$ be reals and $\theta\in ]-\pi/2,\pi/2[-\{0,-\pi/4\}$ be an angle. If $\sigma\in
\BS^{m_1,m_2}_{1,0;\theta}$ then $T_\sigma^{*1}=T_{\sigma_1^*}$ with $\sigma_1^*\in \BS^{(m_1,m_2),1}_{1,0;\theta^{*1}}$ and
$T_\sigma^{*2}=T_{\sigma_2^*}$ with $\sigma_2^*\in \BS^{(m_1,m_2),2}_{1,0;\theta^{*2}}$, where the two angles $\theta^{*1}$
and $\theta^{*2}$ are defined by
$$ \cot(\theta) + \cot(\theta^{*1}) = -1 \qquad \textrm{and} \qquad \tan(\theta)+\tan(\theta^{*2})=-1,$$
and so $\theta^{*1},\theta^{*2} \in ]-\pi/2,\pi/2[-\{0,-\pi/4\}$.
In addition these two symbols $\sigma_1^*$ and $\sigma_1^*$ satisfy an asymptotic formula. We obtain similar results for the two classes $\BS^{m_1,m_2,1}_{1,0;\theta}$ and $\BS^{m_1,m_2,2}_{1,0;\theta}$.
\end{thm}
 
\begin{thm} \label{thm:composition}
 Let $\theta\in]-\pi/2,\pi/2[$ and $\sigma \in \BS^{m_1,m_2}_{1,0;\theta}$ be fixed. Let $\tau_1\in S^{t_1}_{1,0}$ and $\tau_2\in S^{t_2}_{1,0}$ be two linear symbols of order $t_1$ and $t_2$. Then the operator
$$ T(f,g):= T_{\sigma}\left( \tau_1(x,D)f,\tau_2(x,D)g \right) $$
corresponds to the operator $T_m$ with the symbol $m\in \BS_{1,0;\theta}^{m_1+t_1,m_2+t_2}$.
This new symbol satisfies an asymptotic formula. We will give a similar result for the composition ``on the left'' with the other classes $\BS^{m_1,m_2,1}_{1,0;\theta}$ and $\BS^{m_1,m_2,2}_{1,0;\theta}$.
\end{thm}

\gb By following the ideas of \cite{B} and \cite {terwilleger}, it seems possible to obtain identical results for a multidimensionnal problem.

\section{Action of bilinear operators on Sobolev spaces.}

\mb First, we remember our main result of \cite{B}~:

\begin{thm} \label{thmF} Let $\theta$ be a nondegenerate angle : $\theta\in]-\pi/2,\pi/2[ \setminus \{0,-\pi/4\}$.
Let $p,q$ be exponents such that
$$ 0<\frac{1}{r}=\frac{1}{q}+\frac{1}{p} <\frac{3}{2} \textrm{   and   } 1< p,q \leq \infty.$$
Then for all symbol $\sigma\in BS^{0}_{1,0;\theta}$, the operator $T_\sigma$ is continuous from $L^p \times L^q$ into $L^r$.
\end{thm}

\mb This result is a consequence of ``off-diagonal'' estimates for the $x$-independent symbols. The symbolic class  $\BS^{0,0}_{1,0;\theta}$ is a shake between the class $BS^{0}_{1,0;\theta}$ and some Marcinkiewicz condition.
We have to study also first the bilinear Marcinkiewicz multiplier~:

\begin{prop} \label{theocontinuite100} Let $\sigma$ be a bounded function on $\R^2$ such that \be{conditionssing100} \forall b,c\geq 0
\qquad \left|\partial^{b}_\alpha \partial^c_\beta \sigma(\alpha,\beta)\right| \lesssim \left(1+|\alpha|\right)^{-b}\left(1+|\beta|\right)^{-c}. \ee
Let $1<p,q \leq \infty$ exponents such that
$$0<\frac{1}{r}=\frac{1}{q}+\frac{1}{p}.$$
Then for all $\epsilon >0$ the bilinear operator $T_\sigma$ is continuous from $W^{\epsilon,p} \times W^{\epsilon,q}$ to $L^r$.
\end{prop}

\begin{rem} The condition (\ref{conditionssing100}) is related to the degenerate lines $\alpha=0$ and $\beta=0$, that is why we lose some regularity. In addition, as we require some inhomogeneous decay in (\ref{conditionssing100}), it is possible to have some local estimates to describe an ``off-diagonal'' decay. Here we are just interested in global continuities in Sobolev spaces, so we do not give details about this improvement.\\
We know that the $\epsilon>0$ is necessary. In \cite{GK}, L. Grafakos and N. Kalton construct a counter-example for the $\epsilon=0$ result.
\end{rem}

\dem These kind of operators was already studied by L. Grafakos and N. Kalton in \cite{GK}. \\ 
The continuity of $T_\sigma$ from $W^{\epsilon,p} \times W^{\epsilon,q}$ to $L^r$ is equivalent to the continuity of $T_\lambda$ from $L^{p} \times L^{q}$ to $L^r$, with the new symbol
$$\lambda(\alpha,\beta):= \sigma(\alpha,\beta) \left(1+|\alpha|^2\right)^{-\epsilon/2} \left(1+|\beta|^2\right)^{-\epsilon/2}.$$
In fact the new extra weight satisfies~: there exists a constant $c_\epsilon$ such that for all $\alpha,\beta$
$$ \left(1+|\alpha|^2\right)^{-\epsilon/2} \left(1+|\beta|^2\right)^{-\epsilon/2} \leq C_\epsilon \left(1+ \log\left(2+\frac{1+|\beta|}{1+|\alpha|}\right) \right)^{-2}.$$
So it is easy to see that $\lambda$ verifies the following assumptions~: for all $b,c\geq 0$
\be{conditionssing110}  \left|\partial^{b}_\alpha \partial^c_\beta \lambda(\alpha,\beta)\right| \lesssim \left(1+|\alpha|\right)^{-b}\left(1+|\beta|\right)^{-c} \left(1+ \log\left(2+\frac{1+|\beta|}{1+|\alpha|}\right) \right)^{-2} . \ee
Then Theorem 7.4 of \cite{GK} gives us the desired continuity of $T_\lambda$.
\findem

\mb We now have to prove a similar result for $x$-dependent symbols.

\begin{thm} \label{theocontinuite110} Let $\sigma$ be a bounded function on $\R^3$ such that \be{conditionssin100} \forall a,b,c\geq 0
\qquad \left|\partial_x^a\partial^{b}_\alpha \partial^c_\beta \sigma(x,\alpha,\beta)\right| \lesssim \left(1+|\alpha|\right)^{-b}\left(1+|\beta|\right)^{-c}. \ee
Let $1<p,q \leq \infty$ exponents satisfying
$$0<\frac{1}{r}=\frac{1}{q}+\frac{1}{p} <\infty.$$
Then for all $\epsilon >0$ the bilinear operator $T_\sigma$ is continuous from $W^{\epsilon,p} \times W^{\epsilon,q}$ to $L^r$.
\end{thm}

\dem We write in details the proof for $1\leq r$, which allows us to simplify the arguments. We explain the modifications to prove the case $r<1$ in Remark \ref{remark}. \\
The arguments to extend results about $x$-independent symbols to $x$-dependent symbols are nowadays well-known. 
As for the previous Proposition, we have to prove the continuity of $T_\lambda$ from $L^{p} \times L^{q}$ to $L^r$ with 
$$\lambda(x,\alpha,\beta):= \sigma(x,\alpha,\beta) \left(1+|\alpha|^2\right)^{-\epsilon/2} \left(1+|\beta|^2\right)^{-\epsilon/2}.$$
First take a smooth function $\chi$ supported on $[-1,1]$ such that~:
$$ \sum_{i\in\Z} \chi(x-i)^2=1.$$
Then we set
$$ T^i(f,g)(x):= \chi(x-i)^2 T_\lambda(f,g)(x).$$
We can write the previous operator as~:
$$T^i(f,g)(x)= U_{x}^i(f,g)(x),$$
with $U^i$ defined by
\begin{align*}
U_{y}^i(f,g)(x):= \chi(x-i)\chi(y-i)\int_{\R^2} e^{ix(\alpha+\beta)} \widehat{f}(\alpha)\widehat{g}(\beta)\lambda(y,\alpha,\beta) d\alpha d\beta.
\end{align*}
Let $i$ be fixed, and denote $I=[i-1,i+1]$. By using the Sobolev's embedding $W^{1,r}(I) \hookrightarrow L^\infty(I)$ (because $r\geq 1$), we get
$$ \left| T^i(f,g)(x) \right| \leq \| U_{y}(f,g)(x)\|_{\infty,y} \lesssim \sum_{k=0}^{1} \| \partial_y^k U_{y}^i(f,g)(x)\|_{r,I,dy}.$$
Then by integrating for $x\in I$ and using Fubini's Theorem, we obtain
$$ \left\| T^i(f,g)\right\|_{r} \lesssim \sum_{k=0}^{1}  \left\| \left\|\partial_y^k U_{y}^i(f,g) \right\|_{r} \right\|_{r,I,dy}.$$
We can fix $k\in \{0,1\}$ and $y\in [i-1,i+1]$. Then the operator $\partial_y^k U_{y}(f,g)$ corresponds to the operator $T_\tau$ with the $x$-independent symbol
$$\tau(\alpha,\beta):= \chi(y-i)\partial_{y}^{k} \lambda(y,\alpha,\beta).$$
From the assumptions about $\sigma$ (as for Proposition \ref{theocontinuite100}), the symbol $\tau$ satisfies~: for all $b,c\geq 0$
\be{marcin} \left|\partial_\alpha^b \partial_\beta^c \tau(\alpha,\beta) \right|\lesssim \left(1+|\alpha|\right)^{-b}\left(1+|\beta|\right)^{-c}\left(1+ \log\left(2+\frac{1+|\beta|}{1+|\alpha|}\right) \right)^{-2}. \ee
So the operator $T_\tau$ is bounded on Lebesgue spaces. Let now decompose the function $f$ by
$$f:= f{\bf 1}_{4I} + f {\bf 1}_{(4I)^c} := f_0+f_\infty$$ and the function $g$ by
$$g:= g{\bf 1}_{4I} + g {\bf 1}_{(4I)^c} := g_0+g_\infty.$$
With the boundedness of $T_\tau$, we get 
\be{equation1} \left\| T_\tau(f_0,g_0) \right\|_r \lesssim \|f_0\|_{p} \|g_0\|_{q} \lesssim \|f\|_{p,4I} \|g\|_{q,4I}. \ee
We now study the term $T_\tau(f_0,g_\infty)$ and the two other ones can been studied with the same reasoning.
We study the decay of the bilinear kernel of $T_\tau$.
So the kernel $K_\tau(x,y,z)$ is defined by
$$K_\tau(x,y,z):= \int e^{i\alpha(x-y)+i\beta(x-z)} \tau(\alpha,\beta) d\alpha d\beta$$
in a distributionnal sense.
So with integrations by parts, we get that for any integer $M\geq 0$
$$ (x-z)^{2M} K_\tau(x,y,z):= \int e^{i\alpha(x-y)+i\beta(x-z)} \partial_\beta^{2M} \tau(\alpha,\beta) d\alpha d\beta = K_{\partial_\beta^{2M} \tau}(x,y,z).$$
However we have assume some inhomogeneous decay for the symbol $\sigma$ and so for the symbol $\tau$. We also obtain that the symbol $\partial_\beta^{2M}\tau$ satisfies (\ref{marcin}) too. Therefore the operator associated to this symbol satisfies some boundedness on Lebesgue spaces. We write for $x\in I$
\begin{align*}
 \left|T_\tau(f_0,g_\infty)(x) \right| & = \left|\int K_\tau(x,y,z) f_0(y) g_\infty(z) dy dz\right| \\
 & = \left| \int K_{\partial_\beta^{2M} \tau}(x,y,z) f_0(y) (x-z)^{-2M}g_\infty(z) dy dz \right| \\  
 & = \left\| \int K_{\partial_\beta^{2M} \tau}(x,y,z) f_0(y) (w-z)^{-2M}g_\infty(z) dy dz \right\|_{\infty,\omega\in I}. 
\end{align*}
Then by using a Sobolev's embedding and the boundedness of $T_{\partial_\beta^{2M} \tau}$, we get
$$ \left\|T_\tau(f_0,g_\infty)(x) \right\|_{r,I} \lesssim \|f_0\|_{p} \left\|\left\|(w-.)^{-2M}g_\infty\right\|_{q}\right\|_{r,I,dw},$$
with an exponent $N$ as large as we want. However by definition $g_\infty$ is supported on $(4I)^c$, so we can estimate the last norm by
\begin{align*}
  \left\|\left\|(w-.)^{-2M}g_\infty\right\|_{q}\right\|_{r,I,dw} & \lesssim \left\|(1+(i-.)^2)^{-M}g\right\|_{q} \\
 & \lesssim \left(\sum_{k\geq 0} 2^{-kN} \|g\|_{p,2^{k}I}\right),
\end{align*}
where $N$ is an other integer. The first integer $M$ being as large as we want, so the second integer $N$ is too. 
We compute the same arguments for $T_\tau(f_\infty,g_0)$ and $T_\tau(f_\infty,g_\infty)$ and we obtain also
$$ \left\|T^i(f,g)\right\|_{r} \lesssim \left(\sum_{k\geq 0} 2^{-kN} \|f\|_{p,2^{k}I}\right)
\left(\sum_{k\geq 0} 2^{-kN} \|g\|_{q,2^{k}I}\right).$$ 
These we can sum over the index $i$. It is easy to see (as in \cite{B}) that we finally obtain
$$ \left\|T_\lambda(f,g)\right\|_{r} \lesssim \sum_i \left\|T^i(f,g)\right\|_{r} \lesssim \|f\|_p \|g\|_q, $$ 
which is equivalent to the desired continuity.
\findem

\begin{rem} \label{remark} In Proposition \ref{theocontinuite100} for the Marcinkiewicz multipliers, the exponent $r$ may be lower than $1$. In Theorem \ref{theocontinuite110}, we have used a Sobolev's embedding and so the fact that $r\geq 1$. But as we have described in \cite{B}, by using the concept of ``restricted weak type'', we can go around this difficulty and obtain the Theorem \ref{theocontinuite110} for all $r$. The ``restricted weak type'' for the Marcinkiewicz multipliers (associated to $x$-independent symbols) can be proved by the arguments of M. Lacey in \cite{laceyM} or those of C. Muscalu, J. Pipher, T. Tao and C. Thiele in \cite{mptt}. 
\end{rem}

\mb We are now able to have the general result~:

\begin{thm} \label{theocontinuite101} Let $\theta$ be a angle with $\theta\in]-\pi/2,\pi/2[ \setminus \{0,-\pi/4\}$. Let $1<p,q\leq\infty$ be exponents such that
$$0<\frac{1}{r}=\frac{1}{q}+\frac{1}{p} <\frac{3}{2}.$$
Let $\sigma$ be a symbol in the class $\BS^{(0,0)}_{1,0;\theta}$, $\BS^{(0,0),1}_{1,0;\theta}$ or $\BS^{(0,0),2}_{1,0;\theta}$.
Then for all $\epsilon >0$, the operator $T_\sigma$ is continuous from $W^{\epsilon,p} \times W^{\epsilon,q}$ into $L^r$.
\end{thm}

\dem We prove only the first case. The other cases are obtained by the same arguments or by using duality
 (see Subsection \ref{dualite}). In addition for convenience, we deal only with the case $\theta=\pi/4$. \\
 So let $\sigma\in \BS^{(0,0)}_{1,0;\pi/4}$. We use a
smooth partition of the frequency plane in order to restrict the symbol on the three domains $\{ |\alpha|\simeq |\beta|\}$, $\{ |\alpha| <<
|\beta|\}$ and $\{ |\alpha| >> |\beta|\}$. In the first place, the minimums $\min\{|\alpha|,|\alpha-\beta| \}$
and $\min\{|\beta|,|\alpha-\beta|\}$ are equivalent to $|\alpha-\beta|$. So when the symbol is supported in this region,
the result is a consequence of Theorem \ref{thmF} (in this region we can have continuity in Lebesgue spaces : with $\epsilon=0$). In the two other regions, the minimums are respectively equivalent
to $|\alpha|$ and $|\beta|$. So in these two cases, the result is a consequence of Theorem \ref{theocontinuite110}.
\findem

\mb
\begin{rem} In fact, we can show that the operator $T_\sigma$ satisfies some ``off-diagonal'' estimates. By using the definition in
\cite{B}, we can proved that the operator a symbol $\sigma$ of order $(0,0)$, the operator $T_\lambda$ satisfies ``off-diagonal'' estimates at the scale $1$ for any order $\delta\geq 1$ :
$T_\lambda\in {\mathcal O}_{1,\delta}(L^p\times L^q,L^r)$, with the modified symbol
$$\lambda(x,\alpha,\beta):= \sigma(x,\alpha,\beta) \left(1+|\alpha|^2\right)^{-\epsilon/2} \left(1+|\alpha|^2\right)^{-\epsilon/2}.$$
We have also continuities in weighted spaces too, see Theorem 4.4 of \cite{B}.
\end{rem}

\gb We have obtained the main result for bilinear pseudodifferential operators of order $(0,0)$. Before to generalize it for operators of any order, we want to describe other continuity results.

\begin{rem} Let us describe in this remark, an other continuity result using the Modulation spaces. For $0\leq p,t \leq \infty$, man define the Modulation spaces ${\mathcal M}^{p,t}$ (see for example \cite{BO}). We recall the construction of these spaces. Let $\phi\in\s(\R)$ be a non-zero smooth real function. We set the short-time Fourier transform $V_\phi(f)$, given by~: for all $(x,\xi)\in \R^2$
$$V_\phi(f)(x,\xi) := \int_\R e^{-it\xi} f(t) \phi(t-x) dt.$$
Then the space ${\mathcal M}^{p,t}$ is the closure of the Schwartz space for the following norm
$$ \|f\|_{{\mathcal M}^{p,t}}:= \left\| \left\| V_\phi(f)(x,\xi) \right\|_{p,dx} \right\|_{t,d\xi}.$$
It is well known that these spaces do not depend on the function $\phi$. For $p=t=2$, the space ${\mathcal M}^{2,2}$ exactly corresponds to the space $L^2$. \\
With the assumptions of Theorem \ref{theocontinuite101}, the operator $T_\sigma$ is continuous from ${\mathcal M}^{p,t_1} \times {\mathcal M}^{q,t_2}$ into ${\mathcal M}^{r,t_3}$ for all exponents $t_1,t_2,t_3$ satisfying
$$1+\frac{1}{t_3}=\frac{1}{t_1}+\frac{1}{t_2} \quad \textrm{and} \quad 1\leq t_1,t_2,t_3\leq \infty.$$
We also get (for $p=q=t_1=t_2=2$) the continuity from $L^2 \times L^2$ into the Modulation space ${\mathcal M}^{1,\infty}$, which contains the Lebesgue space $L^1$. \\
This result about Modulation spaces comes from Theorem 1 of \cite{BO}, proved by A. Benyi and K. Okoudjou. The use of these spaces allows us to not lose regularity. 
\end{rem}

\begin{rem} To be able to have the $\epsilon=0$ result in Theorem \ref{theocontinuite101}, we have to put an extra assumption about the symbol. We will use a weight (as appeared in \cite{GK})~:
\be{weight} \Omega(\alpha,\beta):=\log\left(2+ \left|\log\left(1+\frac{1+|\beta|}{1+|\alpha|}\right) \right| \right). \ee
For example if we assume that the symbol $\sigma$ satisfies for all $a,b,c\geq 0$
\begin{align*}
\lefteqn{\left| \partial_x^a \partial_\alpha ^b \partial_\beta^c \sigma(x,\alpha,\beta) \right| \lesssim } & & \\
& & \max\left\{\frac{\Omega(\alpha,\beta)}{1+|\alpha|},\frac{1}{1+|\beta-\tan(\theta)\alpha|} \right\}^{b} \max\left\{\frac{\Omega(\alpha,\beta)}{1+|\beta|},\frac{1}{1+|\beta-\tan(\theta)\alpha|}\right\}^{c},
\end{align*}
 then we can take $\epsilon=0$ in Theorem \ref{theocontinuite101}. This result is due to the fact that the additional decay allows us to have $\epsilon=0$ in Proposition \ref{theocontinuite100} (see \cite{GK}).
\end{rem}

\mb Now we can study bilinear operators of any order~:

\begin{thm} \label{thm1} Let $\sigma$ be a bilinear symbol of order $(m_1,m_2)$~: $\sigma\in \BS_{1,0;\theta}^{m_1,m_2}$ with $\theta \in ]-\pi/2,\pi/2[ \setminus\{0,-\pi/4\}$.
Then for all $\epsilon>0$, the bilinear operator $T_\sigma$ is continuous from $W^{\epsilon+(m_1)_+,p} \times W^{\epsilon+(m_2)_+,q}$ in $L^r$ for
all exponents $p,q,r$ satisfying
$$0<\frac{1}{r}=\frac{1}{p}+\frac{1}{q} < \frac{3}{2}$$
and $1< p,q<\infty$.
We have the same result for the classes $\BS_{1,0;\theta}^{(m_1,m_2),1}$ and $\BS_{1,0;\theta}^{(m_1,m_2),2}$.
\end{thm}

\dem \\
As for the previous theorem, we will deal only with the first case and we will assume $\theta=\pi/4$ for convenience. \\
Let $\Phi$ be a smooth function on $\R$ such that
$$ |x|\leq 1 \Longrightarrow \Phi(x)=1 \textrm{  and  supp} (\Phi) \subset B(0,2).$$
We have the frequential decomposition
$$ 1 = \Phi(\alpha-\beta)\Phi(\alpha)\Phi(\beta) + Rest(\alpha,\beta).$$
We set the first symbol
$$ \sigma_1(x,\alpha,\beta):= \sigma(x,\alpha,\beta) \Phi(\alpha-\beta)\Phi(\alpha)\Phi(\beta) $$
and the second one
$$ \sigma_2(x,\alpha,\beta):=\sigma(x,\alpha,\beta) Rest(\alpha,\beta) .$$
$*$ - The study of $\sigma_1$. \\
The function $\Phi$ permits to localize the three quantities $|\alpha|$, $|\beta|$ and $|\alpha-\beta|$. The symbol
$\sigma$ being of order $(m_1,m_2)$, the restriction $\sigma_1$ is also a bilinear symbol of order $(0,0)$ (due to the inhomogeneous assumptions). By applying the previous theorem, we get the desired estimate
$$ \left\| T_{\sigma_1}(f,g) \right\|_{r} \lesssim \|f\|_{W^{\epsilon,p}} \|g\|_{W^{\epsilon,q}} \lesssim \|f\|_{W^{\epsilon+(m_1)_+,p}} \|g\|_{W^{\epsilon+(m_2)_+,q}}.$$
$*$ - The case of $\sigma_2$. \\
The symbol $Rest(\alpha,\beta)$ is composed of 7 terms. Each term is a factor of at least one of these quantities
$1-\Phi(\alpha)$, $1-\Phi(\beta)$ or $1-\Phi(\alpha-\beta)$. All these terms can be studied by the same arguments, so
we deal with the most "extremal"
$$ Rest(\alpha,\beta)= [1-\Phi(\alpha)][1-\Phi(\beta)][1-\Phi(\alpha-\beta)].$$
As $\Phi'$ is supported in a corona around $0$, we have~:
\begin{align}
\left| \partial_\alpha Rest(\alpha,\beta) \right| & \lesssim \left|\Phi'(\alpha)\right| + \left| \Phi'(\alpha-\beta)\right| \nonumber \\
 & \lesssim {\bf 1}_{|\alpha|\simeq 1} + {\bf 1}_{|\alpha-\beta|\simeq 1} \lesssim {\bf 1}_{\min\{|\alpha|,|\alpha-\beta|\} \simeq 1} \nonumber\\
 & \lesssim \left(1+ \min\{|\alpha|,|\alpha-\beta|\} \right)^{-1}. \label{decr1}
\end{align}
We set the symbol
$$ \tau(x,\alpha,\beta) := Rest(\alpha,\beta) \sigma(x,\alpha,\beta) \left(1+ |\alpha|^2 \right)^{-(m_1)_+/2}\left(1+|\beta|^2 \right)^{-(m_2)_+/2}.$$
The symbol $\sigma$ is of order $(m_1,m_2)$, we get also that $\tau$ is a bilinear symbol of order $(0,0)$ with the following estimates (similarly obtained as (\ref{decr1}))~: for all $b,c\geq 0$
\begin{align*}
\left| \partial_\alpha^b \partial_\beta^c Rest(\alpha,\beta) \right| \lesssim  \left(1+ \min\{|\alpha|,|\alpha-\beta|\} \right)^{-b} \left(1+ \min\{|\alpha|,|\alpha-\beta|\} \right)^{-c}.
\end{align*}
We use the new symbol $\tau$, by writing
\begin{align*}
\lefteqn{T_{\sigma_2}(f,g)(x) =} & & \\
 & &  \int_{\R^2} e^{ix(\alpha+\beta)} \widehat{f}(\alpha) \widehat{g}(\beta) \left(1+ |\alpha|^2 \right)^{(m_1)_+/2}\left(1+|\beta|^2 \right)^{(m_2)_+/2} \tau(x,\alpha,\beta) d\alpha d\beta.
\end{align*}
By applying the previous theorem to the symbol $\tau$, we obtain
\begin{align*}
\left\| T_{\sigma_2}(f,g) \right\|_{r} & \lesssim  \left\| (I- \Delta)^{(m_1)_+/2} f\right\|_{W^{\epsilon,p}} \left\| (I- \Delta)^{(m_2)_+/2} g\right\|_{W^{\epsilon,q}} \\ & \lesssim \|f\|_{W^{\epsilon+(m_1)_+,p}} \|g\|_{W^{\epsilon+(m_2)_+,q}}.
\end{align*}
With the same arguments, we can study all the terms due to the decomposition and also we finish the proof.
\findem

\begin{rem} It is interesting to note that we lose some regularity when the orders are negative. We have the same phenomenon with the H\"older's inequality in Sobolev spaces~:
$$  \left\| fg \right \|_{W^{m,r}} \lesssim \|f\|_{W^{m_+,p}} \|g\|_{W^{m_+,q}},$$
which is optimal.
\end{rem}

\mb Now we can prove the following complete result~:

\begin{thm} \label{thm2} Let $\sigma\in \BS_{1,0;\theta}^{m_1,m_2}$ (or $\BS_{1,0;\theta}^{(m_1,m_2),1}$ , $\BS_{1,0;\theta}^{(m_1,m_2),2}$) be a bilinear symbol of order $(m_1,m_2)$ with $\theta \in ]-\pi/2,\pi/2[ \setminus\{0,-\pi/4\}$. Let $1<p,q<\infty$ be exponents satisfying
$$ 0< \frac{1}{r}= \frac{1}{p}+\frac{1}{q} <\frac{3}{2}.$$
Then for all $\epsilon>0$, for all integer $n>0$ and for all functions $f,g\in\s(\R)$~:
\be{derivation} \left\| D^{(n)} T_\sigma(f,g) \right\|_{r} \lesssim \sum_{\genfrac{}{}{0pt}{}{0\leq i,j \leq n}{i+j\leq n}} \|D^{(i)} f\|_{W^{\epsilon+(m_1)_+,p}} \|D^{(j)} g\|_{W^{\epsilon+(m_2)_+,q}}. \ee
Here we set $D^{(i)}$ for the derivation operator of order $i$. Consequently for all $s\geq 0$, we get that $T_\sigma$ is continuous from $W^{s+\epsilon+(m_1)_+,p}\times W^{s+\epsilon+(m_2)_+,q}$ in $W^{s,r}$.
\end{thm}

\dem The link between Theorem \ref{thm1} and Theorem \ref{thm2} is already proved in Proposition 4.7 of \cite{B}. So we do not repeat it and we have also proved Theorem \ref{thm:central}. \findem

\mb By using Sobolev's embedding, there exists an other possible homogeneity in the Lebesgue exponents, which is described in Proposition 3 of \cite{bipseudo}.

\begin{cor} \label{xor} Let $\sigma\in \BS_{1,0;\theta}^{m_1,m_2}$ (or $\BS_{1,0;\theta}^{(m_1,m_2),1}$ , $\BS_{1,0;\theta}^{(m_1,m_2),2}$) be a bilinear symbol of order $(m_1,m_2)$ with $\theta \in ]-\pi/2,\pi/2[ \setminus\{0,-\pi/4\}$. Then for all exponents $1<p,q<\infty$, all $\epsilon>0$ and for all reals $s,t\geq 0$ satisfying
$$ 0<\frac{1}{r_t} = \frac{1}{p} + \frac{1}{q} - t < \frac{3}{2} \quad \textrm{and} \quad 0\leq t\leq \frac{1}{p},\frac{1}{q}\leq 1, $$
the operator $T_\sigma$ is continuous from $W^{s+t+\epsilon+(m_1)_+,p}\times W^{s+t+\epsilon+(m_2)_+,q}$ in $W^{s+t,r_t}$.
\end{cor}

\dem The case $t=0$ corresponds to Theorem \ref{thm2}. The proof for the other cases (as an application of Theorem \ref{thm2}) is written in \cite{bipseudo}.
\findem 

\mb Having obtained a good description of the action of our bilinear operators in Sobolev spaces, we will in the next section study some rules of the symbolic bilinear calculus.

\section{Study of the symbolic calculus for linear and bilinear pseudodifferential operators.}
\label{rules}

In this section, we will study two rules of symbolic calculus. For a function $\sigma:=(x,\xi,\eta) \to \sigma(x,\xi,\eta)\in C^\infty(\R^3)$, we will always denote $\partial_\xi$ and $\partial_\eta$ the two differentiations with respect to the first and the second frequency variable. In addition we write $\partial_{\eta-\xi}$ for the differential operator $\partial_\eta - \partial_\xi$.  

\subsection{The action of duality.}
\label{dualite}

\begin{df} Let $T$ be a bilinear operator acting from $\s(\R) \times \s(\R)$ in $\s'(\R)$. We denote its two adjoints $T^{*1}$ and $T^{*2}$ defined by
$$ \forall f,g,h\in \s(\R^d) \qquad \langle T(f,g),h \rangle = \langle T^{*1}(h,g),f\rangle = \langle T^{*2}(f,h),g\rangle.$$
\end{df}

\mb We recall the duality result for linear pseudodifferential operator (see \cite{AG})~:
\begin{prop} \label{pseudolineaire} Let $\tau\in S^{t}_{1,0}$ be a linear symbol, we define the associated operator
$$\forall f\in\s(\R), \qquad \tau(x,D)(f)(x):= \int_{\R} e^{ix\xi} \widehat{f}(\xi) \tau(x,\xi) d\xi.$$
Then this operator admits an adjoint and $\tau(x,D)^* = \tau^*(x,D)$ with
$$ \tau^{*}(x,\xi)=(2\pi)^{-1}\iint e^{-iy \eta} \overline{\tau}(x-y,\xi-\eta) dyd\eta.$$
We have the asymptotic formula : for all integer $N$,
$$\tau^{*}(x,\xi) - \sum_{k =0}^{N-1} \frac{i^{k}}{k !} \partial_\xi ^k \partial_x^k \overline{\tau(x,\xi)} \in S^{t_1-N}.$$
\end{prop}

\mb In the bilinear case, we have the following result~·

\begin{thm} \label{theodualite} Let $m_1,m_2\in\R$ be reals and $\theta\in ]-\pi/2,\pi/2[-\{0,-\pi/4\}$ be a fixed parameter. If $\sigma\in
\BS^{m_1,m_2}_{1,0;\theta}$ then $T_\sigma^{*1}=T_{\sigma_1^*}$ with $\sigma_1^*\in \BS^{(m_1,m_2),1}_{1,0;\theta^{*1}}$ and
$T_\sigma^{*2}=T_{\sigma_2^*}$ with $\sigma_2^*\in \BS^{(m_1,m_2),2}_{1,0;\theta^{*2}}$, where the two angles $\theta^{*1}$
and $\theta^{*2}$ are defined by
$$ \cot(\theta) + \cot(\theta^{*1}) = -1 \qquad \textrm{and} \qquad \tan(\theta)+\tan(\theta^{*2})=-1.$$
So $\theta^{*1},\theta^{*2} \in ]-\pi/2,\pi/2[-\{0,-\pi/4\}$.
In addition we have
\begin{eqnarray*}
 \sigma_1^*(x,\xi,\eta)  =(2\pi)^{-1}\iint_{\R^2} \overline{\sigma(y,-\alpha-\eta,\eta)} e^{-i(z-\xi)(x-y)} dyd\alpha \\
 \sigma_2^*(x,\xi,\eta)  =(2\pi)^{-1}\iint_{\R^2} \overline{\sigma(y,\xi,-\alpha-\xi)} e^{-i(\alpha-\eta)(x-y)} dyd\alpha.
 \end{eqnarray*}
For these new symbols, we have the asymptotic formulas
$$  \forall N\geq 0, \qquad \sigma_1^*(x,\xi,\eta)-\sum_{k=0}^{N-1} \frac{i^k}{k!} \partial_x^k\partial_\xi^k \overline{\sigma(x,-\xi-\eta,\eta)} \in \BS^{(m_1-N,m_2),1}_{1,0;\theta^{*1}}$$
and
 $$  \forall N\geq 0, \qquad \sigma_2^*(x,\xi,\eta)-\sum_{k=0}^{N-1} \frac{i^k}{k!} \partial_x^k\partial_\eta^k \overline{\sigma(x,\xi,-\eta-\xi)} \in \BS^{(m_1,m_2-N),2}_{1,0;\theta^{*2}}.$$
\end{thm}

\dem In Theorem $4$ of \cite{bipseudo}, it was shown that $T^{*1}$ and $T^{*2}$ are associated to the symbols $\sigma_{1}^*$ and $\sigma_2^*$ given by
\begin{align*}
 \sigma_1^*(x,\xi,\eta) & =(2\pi)^{-1} \iint_{\R^2} \overline{\sigma(y,-\alpha-\eta,\eta)} e^{-i(\alpha-\xi)(x-y)} dyd\alpha \\
 \sigma_2^*(x,\xi,\eta) & =(2\pi)^{-1} \iint_{\R^2} \overline{\sigma(y,\xi,-\alpha-\xi)} e^{-i(\alpha-\eta)(x-y)} dyd\alpha.
\end{align*}
By symmetry, we will only study the first symbol. \\
$1-)$ Proof of $\sigma_1^* \in \BS^{(m_1,m_2),1}_{1,0;\theta^{*1}}$. \\
The symbol $\sigma_1^*$ is given by an oscillating integral, so we recall the main lemma about these integrals (see \cite{AG}).

\begin{lem} \label{oscidualite} For all exponents $m\in\R$, there exists a constant $C=C(m)$ such that
$$ \left| \iint_{\R^2} e^{i\alpha y} a(y,\alpha) dy d\alpha \right| \leq C \left\|a\right\|_{{\mathbb A}^m}$$
with the following norm
$$ \left\|a\right\|_{{\mathbb A}^m} := \sup_{0\leq j,l} \sup_{(y,\alpha)\in\R^2}\left(1+ |y| +|\alpha| \right)^{-m} \left| \partial_y^j \partial_\alpha^l a(y,\alpha) \right|.$$
\end{lem}
\mb By a change of variables, we obtain
$$ \sigma_1^*(x,\xi,\eta)  =\iint_{\R^2} \overline{\sigma(x-y,\alpha-\xi-\eta,\eta)} e^{i\alpha y} dyd\alpha.$$
In addition, it is well known (see \cite{AG}) that we can formally differentiate the oscillating integrals, hence the following estimate
$$ \left| \partial_{x}^{a}\partial_{\xi}^{b}\partial_{\eta-\xi}^{c}  \sigma_1^*(x,\xi,\eta) \right| \lesssim \left\|\partial_{x}^{a}\partial_{\xi}^{b}\partial_{\eta-\xi}^{c} \sigma(x-y,\alpha-\xi-\eta,\eta)\right\|_{(y,\alpha), {\mathbb A}^{|m_1|+|m_2|+b+c}}.$$
By the properties of the symbol $\sigma$, we have
\begin{align*}
 \lefteqn{ \left\| \partial_{x}^{a}\partial_{\xi}^{b}\partial_{\eta-\xi}^{c} \sigma(x-y,\alpha-\xi-\eta,\eta) \right\|_{(y,\alpha),{\mathbb A}^{|m_1|+|m_2|+b+c}} } & & \\
 & \quad  := \sup_{0\leq j,l} \sup_{(y,\alpha)\in\R^2}\left(1+ |y| +|\alpha| \right)^{-(|m_1|+|m_2|+b+c)} \left| \partial_{x}^{a}\partial_{\xi}^{b}\partial_{\eta-\xi}^{c} \partial_y^j \partial_\alpha^l \sigma(x-y,\alpha-\xi-\eta,\eta) \right| \\
 & \quad \lesssim \sup_{0\leq j,l} \sup_{(y,\alpha)\in\R^2}\left(1+ |y| +|\alpha| \right)^{-|m_1|-|m_2|-b-c} \left(1+|\alpha-\xi-\eta|\right)^{m_1}\left(1+|\eta|\right)^{m_2} \\
 & \qquad  \left(1+\min\{|\alpha-\xi-\eta|,|\eta-\tan(\theta)(\alpha-\xi-\eta)| \}\right)^{-b} \\ 
 & \qquad  \left(1+\min\{|\eta|,|\eta-\tan(\theta)(\alpha-\xi-\eta)|\}\right)^{-c} \\
 & \quad \lesssim \left(1+|\xi+\eta|\right)^{m_1}\left(1+|\eta|\right)^{m_2}
 \left(1+\min\{|\xi+\eta|,|\eta-\tan(\theta)(-\xi-\eta)| \}\right)^{-b} \\
 & \qquad \left(1+\min\{|\eta|,|\eta-\tan(\theta)(-\xi-\eta)|\}\right)^{-c} \\
 & \quad \lesssim \left(1+|\xi+\eta|\right)^{m_1}\left(1+|\eta|\right)^{m_2}
\left(1+\min\{|\xi+\eta|,|\eta-\frac{\tan(\theta)}{1+\tan(\theta}\xi| \}\right)^{-b} \\
 & \qquad  \left(1+\min\{|\eta|,|\eta+\frac{\tan(\theta)}{1+\tan(\theta)}\xi)|\}\right)^{-c}.
\end{align*}
Here we have used the Peetre's inequality~:
\be{peetre} \forall s\in\R,\ \forall  u,v\in\R \qquad \left( 1+ |u-v|\right)^{s} \leq \left(1+|u|\right)^{|s|} \left(1+|v|\right)^{s}. \ee
With the definition of $\theta^{*1}$, we also get
\begin{align*}
\lefteqn{ \left| \partial_{x}^{a}\partial_{\xi}^{b}\partial_{\eta-\xi}^{c} \sigma_1^*(x,\xi,\eta) \right| \lesssim  \left(1+|\xi+\eta|\right)^{m_1}\left(1+|\eta|\right)^{m_2} } & & \\
 & & \left(1+\min\{|\xi+\eta|,|\eta-\tan(\theta^{*1})\xi| \}\right)^{-b} \left(1+\min\{|\eta|,|\eta-\tan(\theta^{*1})\xi)|\}\right)^{-c}.
 \end{align*}
This means exactly that $\sigma_1^* \in \BS_{1,0;\theta^{*1}}^{(m_1,m_2),1}$. \\
$2-)$ Asymptotic formula. \\
As for the proof of Proposition \ref{pseudolineaire}, we can show that
$$ \sigma_1^*(x,\xi,\eta)-\sum_{k=0}^{N-1} \frac{i^{k}}{k!} \partial_x^k\partial_\xi^k \overline{\sigma(x,-\xi-\eta,\eta)}$$
corresponds to a bilinear symbol of order $(m_1-N,m_2)$ which has the same properties than the previous one, and so belongs to the desired class. We let the details to the reader.
\findem

\mb From a functional point of view, we have shown that the dual classes of $\BS^{m_1,m_2}_{1,0;\theta}$ are exactly
$\BS^{(m_1,m_2),1}_{1,0;\theta^{*1}}$ and $BS^{(m_1,m_2),2}_{1,0;\theta^{*2}}$. By the same proof, we have the following result~:

\begin{thm} Let $m_1,m_2\in\R$ be reals and $\theta\in ]-\pi/2,\pi/2[-\{0,-\pi/4\}$ be fixed. If $\sigma\in
\BS^{m_1,m_2,1}_{1,\theta}$ then $T_\sigma^{*1}=T_{\sigma_1^*}$ with $\sigma_1^*\in \BS^{m_1,m_2}_{1,0;\theta^{*1}}$ and
$T_\sigma^{*2}=T_{\sigma_2^*}$ with $\sigma_2^*\in \BS^{(m_1,m_2),1}_{1,0;\theta^{*2}}$ where the two angles $\theta^{*1}$ and
$\theta^{*2}$ are defined by
$$ \cot(\theta) + \cot(\theta^{*1}) = -1 \qquad \textrm{and} \qquad  \tan(\theta)+\tan(\theta^{*2})=-1,$$
 and so $\theta^{*1},\theta^{*2} \in ]-\pi/2,\pi/2[-\{0,-\pi/4\}$.
We have the exact formulas
\begin{eqnarray*}
 \sigma_1^*(x,\xi,\eta)  =(2\pi)^{-1}\iint_{\R^2} \overline{\sigma(y,-\alpha-\eta,\eta)} e^{-i(z-\xi)(x-y)} dyd\alpha \\
 \sigma_2^*(x,\xi,\eta)  =(2\pi)^{-1}\iint_{\R^2} \overline{\sigma(y,\xi,-\alpha-\xi)} e^{-i(\alpha-\eta)(x-y)} dyd\alpha.
 \end{eqnarray*}
For each symbol, we have the asymptotic formulas
$$  \forall N\geq 0, \qquad \sigma_1^*(x,\xi,\eta)-\sum_{k=0}^{N-1} \frac{i^k}{k!} \partial_x^k\partial_\xi^k \overline{\sigma(x,-\xi-\eta,\eta)} \in \BS^{m_1-N,m_2}_{1,0;\theta^{*1}}$$
and
 $$  \forall N\geq 0, \qquad \sigma_2^*(x,\xi,\eta)-\sum_{k=0}^{N-1} \frac{i^k}{k!} \partial_x^k\partial_\eta^k \overline{\sigma(x,\xi,-\eta-\xi)} \in \BS^{m_1,m_2-N}_{1,0;\theta^{*2}}.$$
\end{thm}

\mb We can abstract these correspondences by duality
\begin{align*}
 & \left( \BS^{m_1,m_2}_{1,0;\theta} \right)^{*1} = \BS^{(m_1,m_2),1}_{1,0;\theta^{*1}} & \\
 & \left( \BS^{m_1,m_2}_{1,0;\theta} \right)^{*2} = \BS^{(m_1,m_2),2}_{1,0;\theta^{*2}} & \\
 & \left( \BS^{(m_1,m_2),1}_{1,0;\theta} \right)^{*1} = \BS^{m_1,m_2}_{1,0;\theta^{*1}} & \\
 & \left( \BS^{(m_1,m_2),1}_{1,0;\theta} \right)^{*2} = \BS^{(m_1,m_2),1}_{1,0;\theta^{*2}} & \\
 & \left( \BS^{(m_1,m_2),2}_{1,0;\theta} \right)^{*1} = \BS^{(m_1,m_2),2}_{1,0;\theta^{*1}} & \\
 & \left( \BS^{(m_1,m_2),2}_{1,0;\theta} \right)^{*2} = \BS^{m_1,m_2}_{1,0;\theta^{*2}}. &
\end{align*}

\subsection{The action of composition.}

\mb We first remember the rule of composition in the linear case (see \cite{AG})~:

\begin{prop} Let $\tau_1\in S^{t_1}_{1,0}$ and $\tau_2\in S^{t_2}_{1,0}$ be two linear symbols. Then the operator $\tau_1(x,D) \circ \tau_2(x,D)$ corresponds to $\sigma(x,D)$ with the symbol $\sigma\in S^{t_1+t_2}_{1,0}$ given by
$$ \sigma(x,\xi)=(2\pi)^{-1}\iint_{\R^2} e^{-iy \eta} \tau_1(x,\xi-\eta) \tau_2(x-y,\xi) dy d\eta.$$
We have the following asymptotic formula : for all integer $N$,
$$\sigma(x,\xi) - \sum_{\alpha =0}^{N-1} \frac{i^\alpha}{\alpha !} \partial_\xi ^\alpha \tau_1(x,\xi) \partial_x^\alpha \tau_2(x,\xi) \in S^{t_1+t_2-N}_{1,0}$$
\end{prop}

\mb We now come to the main result of this subsection.

\begin{thm} Let $\theta\in ]-\pi/2,\pi/2[$ and $\sigma$ be a symbol of $\BS^{m_1,m_2}_{1,0;\theta}$. Let $\tau_1 \in S^{t_1}_{1,0}$ and $\tau_2 \in S^{t_2}_{1,0}$ be two linear operators of order $t_1$ and $t_2$. Then the operator
$$ T(f,g):= T_\sigma(\tau_1(x,D)f,\tau_2(x,D)g)$$
corresponds to the operator $T_m$ with the symbol $m\in \BS^{m_1+t_1,m_2+t_2}_{1,0;\theta}$ given by
\begin{align*} 
 \lefteqn{m(x,\xi,\eta):=} & & \\
 & & (2\pi)^{-2} \iiiint_{\R^4} e^{i(x-y)(\alpha-\xi)+i(x-z)(\beta-eta)} \sigma(x,\alpha,\beta) \tau_1(y,\xi) \tau_2(z,\eta) d\eta dz d\xi dy.
\end{align*}
This symbol satisfies the following asymptotic formula~: for all $P,N\geq 0$
 $$  m(x,\xi,\eta)-\sum_{k=0}^{N}\sum_{j=0}^P \frac{i^k}{k!} \frac{i^j}{j!} \partial_x^k\tau_1(x,\xi) \partial_x^j\tau_2(x,\eta)  \partial_\xi^k\partial_\eta^j \sigma(x,\xi,\eta) $$ is belonging to the class $\BS^{m_1+t_1-N,m_2+t_2-P}_{1,0;\theta}$.
\end{thm}

\dem Just by convenience, we will only deal with the case $\theta=\pi/4$. By definition, for $l\in\{1,2\}$~:
$$\tau_l(x,D)(f)(x):= \int_{\R} e^{ix\xi} \widehat{f}(\xi) \tau_l(x,\xi) d\xi $$
and
$$ T_\sigma(f,g)(x):= \iint_{\R^2} e^{ix(\alpha+\beta)} \widehat{f}(\alpha) \widehat{g}(\beta) \sigma(x,\alpha,\beta) d\alpha d\beta.$$
It is easy to check that we obtain
\begin{align*}
T(f,g)(x) = (2\pi)^{-2} \iint_{\R^2} e^{ix(\xi+\eta)} m(x,\xi,\eta) \widehat{f}(\xi)\widehat{g}(\eta)\ d\xi d\eta,
\end{align*}
with the new symbol $m$ defined by
\begin{align*} 
\lefteqn{m(x,\xi,\eta):=} & & \\
 & &  (2\pi)^{-2}\iiiint_{\R^4} e^{i(x-y)(\alpha-\xi)+i(x-z)(\beta-\eta)} \sigma(x,\alpha,\beta) \tau_1(y,\xi) \tau_2(z,\eta) d\alpha dz d\beta dy.
\end{align*}
We use a change of variables to write
\begin{align*} 
\lefteqn{m(x,\xi,\eta):=} & & \\
 & & (2\pi)^{-2} \int_{\R^4} e^{iy\alpha+iz\beta} \sigma(x,\xi-\alpha,\eta-\beta) \tau_1(x-y,\xi) \tau_2(x-z,\eta) d\beta dz d\alpha dy.
\end{align*}
We recognize oscillating integrals~:
$$ m(x,\xi,\eta)= (2\pi)^{-1}\int_{\R^2} e^{iy\alpha} \tau_1(x-y,\xi) \lambda_{x,\eta}(\alpha) d y d\alpha,$$
with
$$ \lambda_{x,\eta}(\alpha):= (2\pi)^{-1}\int_{\R^2} e^{iz\beta} \sigma(x,\xi-\alpha,\eta-\beta) \tau_2(x-z,\eta) dz d\beta .$$
{\bf $1-)$ Estimates for $\lambda$.} \\
We have to study $\lambda_{x,\eta}$.
By using Lemma \ref{oscidualite} and by formally differentiating, we get for $a\geq 0$
$$\left|\partial_\alpha^a \partial \lambda_{x,\eta}(\alpha) \right| \lesssim \left\|\partial_\alpha^a \sigma(x,\xi-\alpha,\eta-\beta) \tau_2(x-z,\eta) \right\|_{(z,\beta),{\mathbb A}^{|m_2|+|m_1|+a}}.$$
This norm corresponds to
\begin{align*}
\lefteqn{\left\| \partial_\alpha^a \sigma(x,\xi-\alpha,\eta-\beta) \tau_2(x-z,\eta) \right\|_{(z,\beta),{\mathbb A}^{|m_2|+|m_1|+a}} =} & & \\
 & & \sup_{0\leq j,l} \sup_{(z,\beta)\in\R^2}\left(1+|z|+|\beta| \right)^{-|m_2|-|m_1|-a} \left| \partial_\alpha^a \partial_z^j \partial_\beta^l \sigma(x,\xi-\alpha,\eta-\beta) \tau_2(x-z,\eta) \right|.
 \end{align*}
We can deduce from the properties of $\sigma$ that
\begin{align*}
 \lefteqn{\left| \partial_\alpha^a\partial_\beta^l \sigma(x,\xi-\alpha,\eta-\beta) \right|} & & \\
  & & \lesssim (1+|\xi-\alpha|)^{m_1}(1+|\eta-\beta|)^{m_2}\left(1+\min\{|\xi-\alpha|,|\xi-\alpha-\eta+\beta|\}\right)^{-a} \\
  & & \qquad \qquad \qquad \left(1+\min\{|\eta-\beta|,|\xi-\alpha-\eta+\beta|\}\right)^{-l} \\
  & & \lesssim \left(1+|\xi-\alpha|\right)^{m_1} \left(1+|\eta-\beta|\right)^{m_2} \left(1+\min\{|\xi-\alpha|,|\xi-\alpha-\eta+\beta|\}\right)^{-a},
\end{align*}
and similarly
$$ \left| \partial_z^j \partial_\beta^l \tau_2(x-z,\eta) \right| \lesssim \left(1+ |\eta|\right)^{t_2}.$$
Therefore, we get
\begin{align*}
\lefteqn{\left\| \partial_\alpha^a \sigma(x,\xi-\alpha,\eta-\beta) \tau_2(x-z,\eta) \right\|_{(z,\beta),{\mathbb A}^{|m_2|+|m_1|+a}}} & & \\
 & \qquad \lesssim \sup_{0\leq j,l} \sup_{(z,\beta)\in\R^2}\left(1+ |z| +|\beta| \right)^{-|m_2|-|m_1|-a} \left(1+|\xi-\alpha|\right)^{m_1} \\
 & \qquad  \qquad  \left(1+|\eta-\beta|\right)^{m_2}\left(1+\min\{|\xi-\alpha|,|\xi-\alpha-\eta+\beta|\}\right)^{-a} \left(1+ |\eta|\right)^{t_2} \\
 & \qquad \lesssim \left(1+ |\eta|\right)^{t_2} \left(1+\min\{|\xi-\alpha|,|\xi-\alpha-\eta|\}\right)^{-a} \\
 & \qquad \qquad \left(1+|\xi-\alpha|\right)^{m_1} \left(1+|\eta|\right)^{m_2}.
\end{align*}
We have used Peetre's inequality (\ref{peetre}). Finally, we obtain
\begin{align}
\lefteqn{\left| \partial_\alpha^a \lambda_{x,\eta}(\alpha) \right| \lesssim \left(1+ \min\{|\xi-\alpha|,|\xi-\alpha-\eta|\}\right)^{-a} } & &  \nonumber \\
 & & \left(1+ |\eta|\right)^{t_2} \left(1+|\xi-\alpha|\right)^{m_1} \left(1+|\eta|\right)^{m_2}. \label{est1}
\end{align}
Now, we can estimate the symbol $m$. \\
{\bf $2-)$ Study of $m$.} \\
We recall that $m$ is given by
$$ m(x,\xi,\eta):= \iint_{\R^2} e^{iy\alpha} \tau_1(x-y,\xi) \lambda_{x,\eta}(\alpha) d y d\alpha. $$
By the same arguments using Lemma \ref{oscidualite}, we get
$$\left| m(x,\xi,\eta) \right| \lesssim \left\| \tau_1(x-y,\xi) \lambda_{x,\eta}(\alpha) \right\|_{(y,\alpha),{\mathbb A}^{|m_1|+|m_2|}}.$$
This norm is estimated by
\begin{align*}
\lefteqn{\left\| \tau_1(x-y,\xi) \lambda_{x,\eta}(\alpha) \right\|_{(y,\alpha),{\mathbb A}^{|m_1|+|m_2|}} =} & & \\
 & &  \sup_{0\leq j,l} \sup_{(y,\alpha)\in\R^2}\left(1+ |y| +|\alpha| \right)^{-|m_1|-|m_2|} \left| \partial_y^j \partial_\alpha^l \tau_1(x-y,\xi) \lambda_{x,\eta}(\alpha) \right|.
\end{align*}
By using (\ref{est1}) and Peetre's inequality as previously, we obtain
\begin{align*}
\left\| \tau_1(x-y,\xi) \lambda_{x,\eta}(\alpha) \right\|_{(y,\alpha),{\mathbb A}^{|m_1|}} \lesssim  \left(1+|\xi|\right)^{t_1}\left(1+ |\eta|\right)^{t_2} \left(1+|\xi|\right)^{m_1} \left(1+|\eta|\right)^{m_2}.
\end{align*}
So we have the desired estimate
\begin{align*}
 \left| m(x,\xi,\eta) \right| \lesssim \left(1+|\xi|\right)^{t_1+m_1}\left( 1+|\eta|\right)^{m_2+t_2}.
\end{align*}
We now have to control the action of differentiations on $m$. It is obvious that the spatial differentiation have no importance. For each
differentiation of $m$ with respect to $\xi$, we have an additional weight
$\left(1+\min\{|\xi-\alpha|,|\xi-\alpha-\eta+\beta|\}\right)^{-1}$. The variables $\alpha$ and $\beta$ are controlled
with an appropriate norm $\|\ \|_{\mathbb{A}}$ and so we obtain $\left(1+\min\{|\xi|,|\xi-\eta|\}\right)^{-1}$. We use
the same arguments for the differentiations of $m$ with respect to $\eta$. We finally obtain
\begin{align*}
 \lefteqn{\left| \partial_x^a \partial_\xi^b \partial_\eta^c m(x,\xi,\eta) \right| \lesssim \left(1+|\xi|\right)^{t_1+m_1}\left( 1+|\eta|\right)^{m_2+t_2} } & & \\
 & &  \left(1+\min\{|\xi|,|\xi-\eta|\}\right)^{-b}\left(1+\min\{|\eta|,|\xi-\eta|\}\right)^{-c}.
\end{align*}
This means exactly that $m\in \BS^{m_1+t_1,m_2+t_2}_{1,0,\pi/4}$. \\
{\bf $3-)$ The asymptotic formula.} \\
 For the intermediate symbol $\lambda_{x,\eta}$, by the same arguments as in the linear case, we obtain
$$ \lambda_{x,\eta}(\alpha) \simeq \sum_{k} \frac{i^k}{k!} \partial_{\eta}^k \sigma(x,\xi-\alpha,\eta)  \partial_{x}^k \tau_2(x,\eta). $$
This means that for all integer $N$
$$\gamma_N(\alpha):=\lambda_{x,\eta}(\alpha) - \sum_{k=0}^{N-1} \frac{i^k}{k!} \partial_{\eta}^k \sigma(x,\xi-\alpha,\eta)  \partial_{x}^k \tau_2(x,\eta)$$
satisfies that for all $p\geq 0$
\begin{align*} 
\lefteqn{\left|\partial_{\alpha}^{p} \gamma_N\right| \lesssim \left(1+ |\eta|\right)^{t_2-N} } & & \\
 & & \left(1+ \min\{ |\xi-\alpha|,|\xi-\alpha-\eta|\}\right)^{-p} \left(1+|\xi-\alpha|\right)^{m_1} \left(1+|\eta|\right)^{m_2}.
\end{align*}
In computing the same proof for the other variable, we obtain the bilinear asymptotic formula. We let the details to the reader.
\findem

\mb We have studied the composition from the right of the bilinear operator, by duality we can obtain the following result for the composition from the left~:

\begin{thm} Let $\theta \in ]-\pi/2,\pi/2[$ and $\sigma \in BS^{(m_1,m_2),1}_{1,0;\theta}$ be fixed. Let $\tau\in S^{t}_{1,0}$ be a linear pseudo-differential operator of order $t$. Then the operator
$$ T(f,g):= \tau(x,D) T_{\sigma}\left(f,g \right) $$
corresponds to $T_m$ with the bilinear symbol $m\in \BS^{(m_1+t,m_2),1}_{1,0;\theta}$ given by
$$m(x,\xi,\eta):= (2\pi)^{-1} \int_{\R^2} e^{i(x-y)(\xi+\eta-\alpha)} \sigma(x,\xi,\eta) \tau(x,\alpha) dz d\alpha. $$
This symbol satisfies the asymptotic formula : for all $N\geq 0$
$$ m(x,\xi,\eta) - \sum_{k=0}^{N-1} \frac{i^k}{k!} \partial_{\xi}^{k} \tau(x,\xi+\eta) \partial_{x}^{k} \sigma(x,\xi,\eta) \in \BS^{(m_1+t-N,m_2-N),1}_{1,0;\theta}. $$
We have similar results for the class $\BS^{(m_1,m_2),2}_{1,0;\theta}$.
\end{thm}

\dem  The proof is identical to the previous one. We can deduce this result too by duality with the previous Theorem and Theorem \ref{theodualite}. \findem

\gb
\begin{rem} This result improve Theorem 3 of \cite{bipseudo}, which deals only the particular case $\theta=-\pi/4$.
Here we can have a result for all $\theta$ due to our definition of symbolic classes. Their our larger than the $BS_{1,0;\theta}$, which are considered in \cite{bipseudo}.
\end{rem}

\gb
\begin{rem} The two compositions (from the left and from the right of a bilinear operator) seem to operate differently. We did not succeed in defining a "good" class of bilinear symbol, which would be invariant by these two compositions. \\
There is also a new phenomenom that did not appear in the linear case : the ``principal symbol'' is not invariant by changes of variable. Let us explain what this expression means. Let $\kappa$ be a smooth and nice diffeomorphism of $\R$ (which allows us to change the frequency variable), we have also the $x$ independent symbol $\kappa(\xi)$. Then for all $x$-independent bilinear symbol $\sigma$, we are interested in the operator 
 $$ \widetilde{T_\sigma}(f,g):=\kappa(x,D)^{-1} T_\sigma(\kappa(x,D)(f),\kappa(x,D)(g)).$$
In the theory of the linear pseudodifferential calculus, it is well-know that such a transformation does not change the principal symbol. By using the previous asymptotic expansion, we see that $\widetilde{T_\sigma}=T_m$ with
$$ m(x,\xi,\eta) \simeq \frac{\kappa(\xi)\kappa (\eta)}{\kappa(\xi+\eta)} \sigma(x,\xi,\eta)$$
Here the symbol ``$\simeq$'' has not really a sense because we have not found a symbolic class which is invariant by the two different compositions. However we can note that the new symbol $m$ has not exactly the same behavior than the initial one $\sigma$.
\end{rem}

\end{document}